\documentclass[preprint]{elsarticle}
\newtheorem{Theorem}{Theorem}
\newtheorem{Lemma}[Theorem]{Lemma}
\newtheorem{Example}[Theorem]{Example}

\begin{document}
\title{On the maximum size of minimal definitive quartet sets}
\author{Chris Dowden}
\ead{dowden@lix.polytechnique.fr}
\address{LIX, \'Ecole Polytechnique, 91128 Palaiseau Cedex, France}

\begin{abstract}
In this paper,
we investigate a problem concerning quartets,
which are a particular type of tree on four leaves.
Loosely speaking,
a set of quartets is said to be `definitive' if it completely encapsulates the structure of some larger tree,
and `minimal' if it contains no redundant information.
Here, we address the question of how large a minimal definitive quartet set on $n$ leaves can be,
showing that the maximum size is at least $2n-8$ for all $n \geq 4$.
This is an enjoyable problem to work on,
and we present a pretty construction,
which employs symmetry.
\end{abstract}


\maketitle

\section{Introduction}

The motivation for this paper comes from the field of phylogenetics,
which involves the study of the `tree of life' depicting all living things,
as popularised by Charles Darwin.
In such a representation,
existing species are drawn as leaves of the tree,
while their ancestors are shown as interior vertices.

In practice,
the overall evolutionary (or `phylogenetic') tree is built up
by piecing together various smaller items of information.
For example,
if species $u$ and $v$ both have wings and species $w$ and $x$ do not,
then it is likely that $u$ and $v$ have a common ancestor that is not shared by $w$ and $x$,
and so the path from $u$ to $v$ on the tree of life should not intersect the path from $w$ to $x$.

The objective of this paper is to present a new result on quartets,
which are a type of graph often used when reconstructing evolutionary trees in this way.
We start by providing some necessary definitions.

\begin{figure} [ht]
\setlength{\unitlength}{1cm}
\begin{picture}(10,2.4)(-6,-0.95)
\put(0,0){\line(0,1){0.5}}
\put(0,0){\line(1,-1){0.7}}
\put(0,0){\line(-1,-1){0.7}}
\put(0,0.5){\line(1,1){0.7}}
\put(0,0.5){\line(-1,1){0.7}}
\put(-0.95,1.25){$u$}
\put(0.8,1.25){$v$}
\put(-0.95,-0.95){$w$}
\put(0.8,-0.95){$x$}
\end{picture}
\caption{The quartet $uv|wx$.} \label{quartet}
\end{figure}
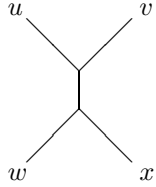

A \emph{phylogenetic tree} is a tree with no vertices of degree $2$
in which the leaves are labelled (distinctly) and the interior vertices are not.
A phylogenetic tree is called \emph{binary} if all interior vertices have degree exactly $3$,
and a \emph{quartet} is defined to be a binary phylogenetic tree with precisely four leaves
(note that such a tree is unique up to labelling).
We use the notation $uv|wx$ to denote a quartet that is labelled as in Figure~\ref{quartet}.

We say that a phylogenetic tree $T$ \emph{displays} the quartet $uv|wx$
if $u$, $v$, $w$ and $x$ are all leaves in $T$
and the path from $u$ to $v$ does not intersect the path from $w$ to $x$
(or, equivalently, there exists a cut-edge in $T$ which separates $u$ and $v$ from $w$ and $x$).
We say that $T$ displays a set of quartets $Q$ if $T$ displays each individual quartet $q \in Q$.

\begin{figure} [ht]

\setlength{\unitlength}{1cm}
\begin{picture}(10,2.4)(-1.5,-0.95)
\put(0,0){\line(0,1){0.5}}
\put(0,0){\line(1,-1){0.7}}
\put(0,0){\line(-1,-1){0.7}}
\put(0,0.5){\line(1,1){0.7}}
\put(0,0.5){\line(-1,1){0.7}}
\put(-0.35,-0.35){\line(1,-1){0.4}}
\put(-0.95,1.25){$u$}
\put(0.8,1.25){$v$}
\put(-0.95,-0.95){$w$}
\put(0.75,-0.95){$x$}
\put(0.1,-0.95){$y$}

\put(3,0){\line(0,1){0.5}}
\put(3,0){\line(1,-1){0.7}}
\put(3,0){\line(-1,-1){0.7}}
\put(3,0.5){\line(1,1){0.7}}
\put(3,0.5){\line(-1,1){0.7}}
\put(3.35,-0.35){\line(-1,-1){0.4}}
\put(2.05,1.25){$u$}
\put(3.8,1.25){$v$}
\put(2.05,-0.95){$w$}
\put(3.75,-0.95){$x$}
\put(2.8,-0.95){$y$}

\put(6,-0.7){\line(0,1){1.2}}
\put(6,0){\line(1,-1){0.7}}
\put(6,0){\line(-1,-1){0.7}}
\put(6,0.5){\line(1,1){0.7}}
\put(6,0.5){\line(-1,1){0.7}}
\put(5.055,1.25){$u$}
\put(6.8,1.25){$v$}
\put(5.05,-0.95){$w$}
\put(6.75,-0.95){$x$}
\put(5.9,-0.95){$y$}

\put(9,0){\line(0,1){0.5}}
\put(9,0){\line(1,-1){0.7}}
\put(9,0){\line(-1,-1){0.7}}
\put(9,0.5){\line(1,1){0.7}}
\put(9,0.5){\line(-1,1){0.7}}
\put(9,0.25){\line(1,0){0.7}}
\put(8.05,1.25){$u$}
\put(9.8,1.25){$v$}
\put(8.05,-0.95){$w$}
\put(9.75,-0.95){$x$}
\put(9.8,0.15){$y$}
\end{picture}
\caption{Four different phylogenetic trees displaying the quartets $uv|wx$ and $uv|wy$.} \label{define}
\end{figure}
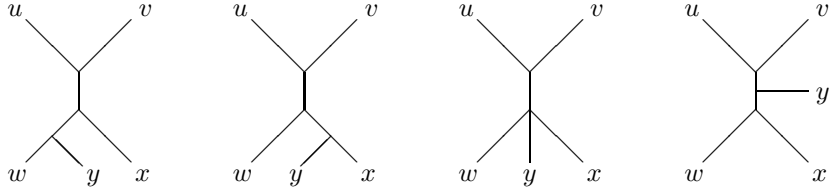

For a set of quartets $Q$ with total leaf-set $L(Q)$,
we say that $Q$ \emph{defines} a tree $T$
(or that $Q$ is \emph{definitive} for $T$)
if $T$ is the unique phylogenetic tree with leaf-set $L(Q)$ that displays $Q$.
Note that many quartet sets will not define any tree,
either because they contain quartets that are incompatible
(e.g.~$\{uv|wx, uw|vx\}$)
or because they are not informative enough to be particular to one tree
(e.g.~$\{uv|wx, uv|wy\}$ is displayed by four different phylogenetic trees with leaf-set $\{u,v,w,x,y\}$,
as shown in Figure~\ref{define}).
An example of a quartet set that is definitive is $\{uv|wx, ux|wy\}$,
which can be seen to define the left-most tree in Figure~\ref{define}.
Finally, we say that $Q$ is a \emph{minimal} definitive quartet set
(or that $Q$ is \emph{minimally definitive})
if $Q$ defines some tree $T$ but,
for all $q \in Q$,
$Q - q$ does not define $T$
(for example, $\{uv|wx, ux|wy\}$ is minimally definitive,
but not $\{uv|wx, ux|wy, uv|wy\}$).

It is fairly straightforward to see that if $Q$ defines $T$,
then $T$ must be binary and $Q$ must \emph{distinguish} every interior edge of $T$
(a quartet $uv|wx$ is said to distinguish an edge $e \in T$ if $e$ is the unique cut-edge in $T$
that separates $u$ and $v$ from $w$ and $x$).
Thus, since a binary tree on $n$ leaves always has exactly $n-3$ interior edges,
it follows that every definitive quartet set on $\{1,2, \ldots, n\}$
must contain at least $n-3$ quartets.
Furthermore, it is known that for any binary phylogenetic tree $T$ on $n$ leaves,
there is indeed a set of $n-3$ quartets that does define $T$
(see, for example, \cite{sem} Corollary 6.3.10).

Hence, the remaining interest in minimal definitive quartet sets lies in the question of how large they can be.
Examples have been produced that have size greater than $n-3$,
but until recently it was thought that the maximum possible size would be bounded by $n+c$ for some fixed constant $c$.
However, Humphries (\cite{hum}, Theorem 3.4.1)
then proved that there actually exist examples with size at least $\frac{3}{2}n - 6$, for all $n \geq 4$.
In this paper,
we will improve matters still further by constructing minimal definitive quartet sets of size $2n-8$, for all $n \geq 5$.

\section{Main Section}

This section will culminate in the inductive construction of minimal definitive quartet sets of size $2n-8$.
The structure of the section will be as follows:
we shall start by stating three lemmas that will be useful to us;
we shall then prove the result for $n=6$,
which will be the base case for our induction;
we shall then also prove the $n=7$ case,
as a way to convey the ideas of the inductive step;
and finally we shall prove the full result.

We start by making explicit a result that we have already noted:

\begin{Lemma}[\cite{sem}, Proposition 6.8.4] \label{interior}
Let $Q$ be a set of quartets that defines a tree $T$.
Then each interior edge of $T$ must be distinguished by at least one quartet in $Q$.
\end{Lemma}

The converse of Lemma~\ref{interior} is known not to be true in general.
However, the following result, which will play an extremely important role in our construction, comes close:

\begin{Lemma}[\cite{sem}, Theorem 6.8.8] \label{intersection}
Let $Q$ be a set of quartets containing a common leaf,
and let $T$ be a tree displaying $Q$ for which each interior edge is distinguished by at least one quartet in $Q$.
Then $Q$ defines $T$.
\end{Lemma}

Often, a set of quartets $Q$ can be used to `infer' a further quartet $q$,
in the sense that every phylogenetic tree displaying $Q$ must also display $q$.
The notation $Q \vdash q$ is used to denote such inferences.
Numerous examples are known,
but we shall only use one very simple one:

\begin{Lemma} \label{inference}
$\{ab|de, bc|de\} \vdash ac|de$.
\end{Lemma}

\begin{figure} [ht]
\setlength{\unitlength}{1cm}
\begin{picture}(10,1.9)(-3.375,-1.1)
\put(0,0){\line(1,0){2.25}}
\put(2.5,0){\line(1,0){0.25}}
\put(3,0){\line(1,0){0.25}}
\put(3.5,0){\line(1,0){0.25}}
\put(4,0){\line(1,0){1.25}}
\put(1,0){\line(0,-1){0.7}}
\put(2,0){\line(0,-1){0.7}}
\put(4.25,0){\line(0,-1){0.7}}
\put(0,0){\line(-1,1){0.6}}
\put(0,0){\line(-1,-1){0.6}}
\put(5.25,0){\line(1,1){0.6}}
\put(5.25,0){\line(1,-1){0.6}}
\put(-0.9,0.55){$1$}
\put(-0.9,-0.75){$2$}
\put(6.05,0.55){$i$}
\put(6.05,-0.75){$i-1$}
\put(0.925,-1.1){$3$}
\put(1.925,-1.1){$4$}
\put(3.9,-1.1){$i-2$}
\end{picture}
\caption{The caterpillar tree $T_{i}$.} \label{Ti}
\end{figure}
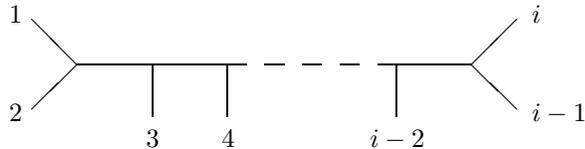

As well as the three lemmas that we have stated,
\emph{caterpillar} trees will also play a major role in our proofs.
The caterpillar tree on $i$ leaves, which we shall denote by $T_{i}$, is defined via Figure~\ref{Ti}.

We shall now give an example of a minimal definitive quartet set on six leaves that has size four,
thus fulfilling our $2n-8$ target.
Such sets have already been produced before now,
but ours has a nice reversible symmetry to it that will later prove significant.

\begin{Lemma} \label{n=6}
The set of quartets
\begin{displaymath}
Q_{6} = \{ 12|35, 13|46, 12|56, 24|56 \}
\end{displaymath}
is minimally definitive for the caterpillar tree $T_{6}$.
\end{Lemma}
\textbf{Proof}
Let us first show that $Q_{6}$ is definitive.
Note that we have $\{12|56, 24|56 \} \vdash 14|56$,
by Lemma~\ref{inference},
and hence $Q_{6} \vdash \{12|35, 13|46, 14|56\}$.
Since the three quartets in this subset all contain the common leaf $1$
and collectively distinguish each interior edge of $T_{6}$,
definitiveness then follows automatically from Lemma~\ref{intersection}.

\begin{figure} [ht]
\setlength{\unitlength}{1cm}
\begin{picture}(10,1.9)(-4.5,-1.1)
\put(0,0){\line(1,0){3}}
\put(1,0){\line(0,-1){0.7}}
\put(2,0){\line(0,-1){0.7}}
\put(0,0){\line(-1,1){0.6}}
\put(0,0){\line(-1,-1){0.6}}
\put(3,0){\line(1,1){0.6}}
\put(3,0){\line(1,-1){0.6}}
\put(-0.9,0.55){$2$}
\put(-0.9,-0.75){$4$}
\put(3.8,0.55){$5$}
\put(3.8,-0.75){$3$}
\put(0.925,-1.1){$6$}
\put(1.925,-1.1){$1$}
\end{picture}
\caption{A tree $T^{\prime}$ displaying the quartet set $Q_{6}-12|56$.} \label{Tprime}
\end{figure}
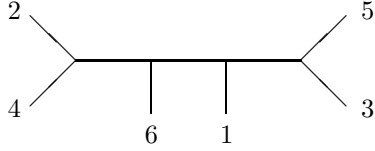

It remains to show that $Q_{6}$ is \emph{minimally} definitive.
If not, then there exists a quartet $q \in Q_{6}$ such that $Q_{6}-q$ defines $T_{6}$.
By Lemma~\ref{interior},
the only possibility for $q$ is $12|56$.
However, the tree $T^{\prime}$ shown in Figure~\ref{Tprime} displays $Q_{6} - 12|56$,
and $T^{\prime}$ is certainly distinct from $T_{6}$.
Hence, it follows that $Q_{6}$ is indeed minimally definitive.
\phantom{qwerty}
\setlength{\unitlength}{.25cm}
\begin{picture}(1,1)
\put(0,0){\line(1,0){1}}
\put(0,0){\line(0,1){1}}
\put(1,1){\line(-1,0){1}}
\put(1,1){\line(0,-1){1}}
\end{picture} \\

We shall now see how to use the set $Q_{6}$ from Lemma~\ref{n=6}
to produce a minimal definitive quartet set on seven leaves that has size six.
This example, combined with the paragraph of discussion after the proof,
is intended to help make clear the general strategy,
which utilises the symmetry properties that we have observed,
but the reader is free to proceed straight to the full proof of Theorem~\ref{2n-8} if he so wishes.

\begin{Example} \label{n=7}
The set of quartets
\begin{displaymath}
Q_{7} = \{ 12|35, 13|46, 12|57, 24|57, 13|67, 35|67 \}
\end{displaymath}
is minimally definitive for the caterpillar tree $T_{7}$.
\end{Example}
\textbf{Proof}
As a rigorous proof is implicitly included within that to Theorem~\ref{2n-8},
we shall provide a slightly more informal treatment here.
Firstly, note that Table~\ref{table}
\begin{table} [ht]
\setlength{\unitlength}{1cm}
\begin{picture}(10,2.2)(-4.5,-0.95)
\begin{tabular}{l}
$\left. \begin{array}{l} 12|35 \\
13|46 \end{array} \right.$ \\
$\left. \begin{array}{l} 12|57 \\
24|57 \end{array} \right\} \vdash 14|57$ \\
$\left. \begin{array}{l} 13|67 \\
35|67 \end{array} \right\} \vdash 15|67$ \\
\end{tabular}
\end{picture}
\caption{Some inferences that can be made from the quartet set $Q_{7}$.} \label{table}
\end{table}
shows that we can again use Lemma~\ref{intersection} to prove definitiveness
(it is worth observing the way that the quartets are paired up here).
By Lemma~\ref{interior},
it then only remains to provide suitable trees $T^{\prime\prime}$ and $T^{\prime\prime\prime}$
displaying $Q_{7} - 12|57$ and $Q_{7} - 13|67$, respectively.
But note that $T^{\prime\prime}$ (see Figure~\ref{n7fig})
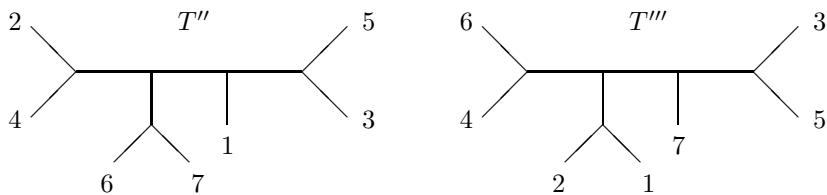
\begin{figure} [ht]
\setlength{\unitlength}{1cm}
\begin{picture}(10,2.4)(-1.4,-1.6)
\put(0,0){\line(1,0){3}}
\put(1,0){\line(0,-1){0.7}}
\put(2,0){\line(0,-1){0.7}}
\put(0,0){\line(-1,1){0.6}}
\put(0,0){\line(-1,-1){0.6}}
\put(3,0){\line(1,1){0.6}}
\put(3,0){\line(1,-1){0.6}}
\put(-0.9,0.55){$2$}
\put(-0.9,-0.75){$4$}
\put(3.8,0.55){$5$}
\put(3.8,-0.75){$3$}
\put(0.325,-1.6){$6$}
\put(1.525,-1.6){$7$}
\put(1.925,-1.1){$1$}
\put(1.35,0.55){$T^{\prime\prime}$}
\put(1,-0.7){\line(-1,-1){0.5}}
\put(1,-0.7){\line(1,-1){0.5}}

\put(6,0){\line(1,0){3}}
\put(7,0){\line(0,-1){0.7}}
\put(8,0){\line(0,-1){0.7}}
\put(6,0){\line(-1,1){0.6}}
\put(6,0){\line(-1,-1){0.6}}
\put(9,0){\line(1,1){0.6}}
\put(9,0){\line(1,-1){0.6}}
\put(5.1,0.55){$6$}
\put(5.1,-0.75){$4$}
\put(9.8,0.55){$3$}
\put(9.8,-0.75){$5$}
\put(6.325,-1.6){$2$}
\put(7.525,-1.6){$1$}
\put(7.925,-1.1){$7$}
\put(7,-0.7){\line(-1,-1){0.5}}
\put(7,-0.7){\line(1,-1){0.5}}
\put(7.35,0.55){$T^{\prime\prime\prime}$}
\end{picture}
\caption{Trees $T^{\prime\prime}$ and $T^{\prime\prime\prime}$
displaying the quartets $Q_{7} - 12|57$ and $Q_{7} - 13|67$, respectively.} \label{n7fig}
\end{figure}
can easily be formed from the tree $T^{\prime}$ in Figure~\ref{Tprime}
(it is important to note the similarity between $Q_{6}$ and $Q_{7}$),
while the symmetry of $Q_{7}$ allows us to take $T^{\prime\prime\prime}$
to be the same as $T^{\prime\prime}$,
but with the numbers reversed!
\phantom{qwerty}
\setlength{\unitlength}{.25cm}
\begin{picture}(1,1)
\put(0,0){\line(1,0){1}}
\put(0,0){\line(0,1){1}}
\put(1,1){\line(-1,0){1}}
\put(1,1){\line(0,-1){1}}
\end{picture} \\

As we shall now see,
the inductive step in the full proof follows exactly the same procedure as in the example above.
We shall always use caterpillars,
and we shall show definitiveness by always adding an extra pair of quartets
that together infer $1(n-2)|(n-1)n$.
Proving minimality will then come down to constructing various trees of the form $Q-q$.
All but one of these will be formed from the trees of the previous stage of the induction,
while the additional tree will be created by reversing numbers.

\begin{Theorem} \label{2n-8}
Let $n \geq 5$ be some positive integer.
Then there exists a minimal definitive quartet set on $n$ leaves that has size $2n-8$.
\end{Theorem}
\textbf{Proof}
We have already noted the result for $n \in \{5,6\}$
(and also for $n=7$, for those that have read Example~\ref{n=7}),
so we shall now proceed inductively,
using the set $Q_{6}$ defined in Lemma~\ref{n=6} as our base.
Let us use $q_{6,1}$, $q_{6,2}$, $q_{6,3}$ and $q_{6,4}$, respectively,
to denote the quartets $12|35$, $13|46$, $12|56$ and $24|56$, in that order
(so $Q_{6} = \{q_{6,1}, q_{6,2}, q_{6,3}, q_{6,4}\}$).
For $k \geq 7$,
let us then define $Q_{k} = \{q_{k,1}, q_{k,2}, \ldots, q_{k,2k-8}\}$
recursively from $Q_{k-1}$ as follows:
(i)~for all $i \leq 2k-12$, set $q_{k,i} = q_{k-1,i}$;
(ii)~for $i \in \{2k-11,2k-10\}$,
set $q_{k,i}$ to be the same as $q_{k-1,i}$,
but with occurrences of $k-1$ replaced by $k$;
and (iii)~set $q_{k,2k-9} = 1(k-4)|(k-1)k$
and $q_{k,2k-8} = (k-4)(k-2)|(k-1)k$.
It can be checked that this procedure produces the set $Q_{7}$ defined in Example~\ref{n=7}.

Note that $|Q_{k}| = 2k-8$ and so, by induction,
it now suffices to prove that $Q_{k}$ is minimally definitive for the caterpillar tree $T_{k}$
given that $Q_{k-1}$ is minimally definitive for $T_{k-1}$.
This is precisely what we shall now do.

First, let us check that $T_{k}$ displays $Q_{k}$.
Note that $T_{k}$ displays all quartets $uv|wx$ for $u<v<w<x$,
and $q_{k,2k-9}$ and $q_{k,2k-8}$ are certainly of this form.
By induction,
we can see that all other quartets in $Q_{k}$ also satisfy this property,
and so $T_{k}$ does indeed display $Q_{k}$.

Next, we shall show that $Q_{k}$ is definitive for $T_{k}$,
using the same argument as for when $k=6$.
By induction, we may assume that
$\{q_{k,1},q_{k,2}, \ldots, q_{k,2k-12}\} \vdash
\{ 12|35,13|46,14|57, \ldots, 1(k-4)|(k-3)(k-1)\}$.
Note $q_{k,2k-11} = 1(k-5)|(k-2)k$
and $q_{k,2k-10} = (k-5)(k-3)|(k-2)k$,
so $\{q_{k,2k-11},q_{k,2k-10}\} \vdash 1(k-3)|(k-2)k$
by Lemma~\ref{inference} (and so the induction does hold).
Finally, Lemma~\ref{inference} also implies that
$\{q_{k,2k-9},q_{k,2k-8}\} \vdash 1(k-2)|(k-1)k$.
Hence, just as with the case when $k=6$,
we may use Lemma~\ref{intersection} to deduce that $Q_{k}$ is definitive for $T_{k}$.

It now only remains for us to show that $Q_{k}$ is \emph{minimally} definitive.
To do this, we need to show that for each $q_{k,i}$ there exists a tree $T_{k,i} \neq T_{k}$ that displays $Q_{k} - q_{k,i}$.

For $i \leq 2k-10$,
we can take $T_{k,i}$ to be the tree formed from $T_{k-1,i}$ by replacing vertex $k-1$ with a `cherry' $\{k-1,k\}$
(by which we mean a rooted binary tree with leaf-set $\{k-1,k\}$ ---
for example, the tree $T^{\prime\prime}$ in Figure~\ref{n7fig}
is formed from the tree $T^{\prime}$ in Figure~\ref{Tprime}
by replacing vertex $6$ with the cherry $\{6,7\}$).
It is clear that $T_{k,i} \neq T_{k}$, since $T_{k-1,i} \neq T_{k-1}$.
The proof that $T_{k,i}$ displays $Q_{k} - q_{k,i}$ follows from observing that
(a)~$T_{k,i}$ displays all quartets that are displayed by $T_{k-1,i}$,
since $T_{k-1,i}$ is a subgraph of $T_{k,i}$,
(b)~for $w<k-1$,
$T_{k,i}$ displays the quartet $uv|wk$ if it displays $uv|w(k-1)$
(since $\{k-1,k\}$ forms a cherry),
and hence $T_{k,i}$ displays $uv|wk$ if $T_{k-1,i}$ displays $uv|w(k-1)$,
and (c)~$T_{k,i}$ displays all quartets of the form $uv|(k-1)k$,
again using the fact that $\{k-1,k\}$ forms a cherry.

For $i=2k-8$,
we may appeal to Lemma~\ref{interior},
and so this only leaves the case $i=2k-9$,
for which we have the quartet
$q_{k,2k-9} = 1(k-4)|(k-1)k$.
To deal with this,
let us take $T_{k,2k-9}$ to be the tree formed from $T_{k,3}$ (which displays $Q_{k} - 12|57$)
by `reversing' all the numbers,
i.e.~$1$ becomes $k$, $2$ becomes $k-1$, $3$ becomes $k-2$, and so on.
Note that $T_{k,2k-9} \neq T_{k}$,
since $T_{k}$ is the `reverse' of itself,
and so it only remains to show that $T_{k,2k-9}$ displays $Q_{k} - q_{k,2k-9}$,
which we shall now do.

First, note that the tree $T_{k,3}$ was formed by taking a tree displaying $Q_{6} - 12|56$,
replacing vertex $6$ with a cherry $\{6,7\}$,
then replacing vertex $7$ with a cherry $\{7,8\}$,
replacing vertex $8$ with a cherry $\{8,9\}$,
and so on until replacing vertex $k-1$ with a cherry $\{k-1,k\}$.
Hence, $T_{k,3}$ must display every quartet of the form $uv|wx$ for $u<v<w<x$ and $w \geq 6$,
and so $T_{k,2k-9}$ must display every quartet of the form $ab|cd$ for $a<b<c<d$ and $b \leq k-5$.

This immediately covers every quartet in $Q_{k} - q_{k,2k-9}$ apart from three:
$q_{k,2k-12} = (k-6)(k-4)|(k-3)(k-1)$,
$q_{k,2k-10} = (k-5)(k-3)|(k-2)k$
and $q_{k,2k-8} = (k-4)(k-2)|(k-1)k$.
Furthermore, since these are the `opposites' of
$q_{k,4} = 24|57$, $q_{k,2} = 13|46$ and $q_{k,1} = 12|35$,
which are all displayed by $T_{k,3}$,
we find that these three remaining quartets are also all displayed by $T_{k,2k-9}$.
Hence, $T_{k,2k-9}$ displays $Q_{k} - q_{k,2k-9}$, and so we are done.
\phantom{qwerty}
\setlength{\unitlength}{.25cm}
\begin{picture}(1,1)
\put(0,0){\line(1,0){1}}
\put(0,0){\line(0,1){1}}
\put(1,1){\line(-1,0){1}}
\put(1,1){\line(0,-1){1}}
\end{picture}

\section{Questions}

The obvious question is whether $2n-8$ can be improved upon,
and it would be interesting to know of any better examples.
Throughout this paper,
we have only used caterpillar trees,
partly for simplicity,
and so another question of interest would be to ask whether caterpillars can always be relied upon
to provide the extremal cases.
Finally, it would also be nice to obtain some sort of upper bound
on the maximum possible size of a minimal definitive quartet set,
other than the trivial $\left( ^{n} _{4} \right)$.

\section*{Acknowledgements}

I would like to thank Charles Semple for introducing me to the problem,
and the reviewers for their helpful comments.

\end{document}